\documentclass[3p]{elsarticle}
\usepackage[all]{xy}
\usepackage{color}
\usepackage{url}
\usepackage{amscd, amssymb, amsfonts, amsthm, amsmath}
\usepackage{mathdots}
\usepackage{hyperref}
\usepackage{multicol}

\newcommand{\soc}{\mbox{\rm{Soc}}}
\newcommand{\topp}{\mbox{\rm{Top}}}
\newcommand{\length}{\mbox{\rm{l}}}
\newcommand{\loewy}{\mbox{\rm{ll}}}
\newcommand{\gl}{\mbox{\rm{gldim}}}
\newcommand{\fin}{\mbox{\rm{findim}}}
\newcommand{\Fin}{\mbox{\rm{Findim}}}

\newcommand{\fidim}{\mbox{\rm{$\phi$dim}}}

\newcommand{\findim}{\mbox{\rm{findim}}}

\def\mod{\mbox{\rm{mod}}}
\def\Mod{\mbox{\rm{Mod}}}
\def\ind{\mbox{\rm{ind}}}
\def\add{\mbox{\rm{add}}}

\def\pd{\mbox{\rm{pd}}}

\def\enn{\hbox{\rm{End}}}
\def\rk{\hbox{\rm{rk}}}

\def\Orb{\mbox{\rm{Orb}}}

\def\Rep{\mbox{\rm{Rep}}}
\def\rep{\mbox{\rm{rep}}}

\def\start{\mbox{\rm{s}}}
\def\target{\mbox{\rm{t}}}

\makeatletter
\def\ps@pprintTitle{%
  \let\@oddhead\@empty
  \let\@evenhead\@empty
  \let\@oddfoot\@empty
  \let\@evenfoot\@oddfoot
}
\makeatother

\begin{document}
\newcommand{\mono}[1]{%
\gdef\puA{#1}}
\newcommand{\puA}{}
\newcommand{\faculty}[1]{%
\gdef\puC{#1}}
\newcommand{\puC}{}
\newcommand{\facultad}[1]{%
\gdef\puD{#1}}
\newcommand{\puD}{}
\newcommand{\N}{\mathbb{N}}
\newcommand{\Z}{\mathbb{Z}}
\newtheorem{teo}{Theorem}[section]
\newtheorem{prop}[teo]{Proposition}
\newtheorem{lema}[teo] {Lemma}
\newdefinition{ej}[teo]{Example}
\newtheorem{obs}[teo]{Remark}
\newtheorem{defi}[teo]{Definition}
\newtheorem{coro}[teo]{Corollary}
\newtheorem{nota}[teo]{Notation}



\title{Igusa-Todorov and LIT algebras on Morita context algebras}

\author[add]{Marcos Barrios}
\ead{marcosb@fing.edu.uy} 

\author[add]{Gustavo Mata\corref{cor}}
\ead{gmata@fing.edu.uy}

\address[add]{Universidad de La Rep\'ublica, Facultad de Ingenier\'ia -  Av. Julio Herrera y Reissig 565, Montevideo, Uruguay}
\cortext[cor]{Corresponding Author}

\begin{abstract}
In this article, we prove that, under certain conditions, Morita context algebras that arise from Igusa-Todorov (LIT) algebras and have zero bimodule morphisms are also Igusa-Todorov (LIT). For a finite dimensional algebra $A$, we prove that the class $\phi_0^{-1}(A) = \{M: \phi(M)=0\}$ is a 0-Igusa-Todorov subcategory if and only if $A$ is selfinjective or $\gl(A)< \infty$. As a consequence $A$ is an $(n,V, \phi_0^{-1}(A))$ algebra if and only if $A$ is selfinjective or $\gl(A)< \infty$. We also show that the opposite algebra of a LIT algebra is not LIT in general. 
\end{abstract} 


\maketitle

\section{Introduction}

For an Artin algebra $A$, are defined the finitistic dimensions as follows

\begin{itemize}
\item $\findim(A) = \sup\{ \pd(M) |\ M \in \mod A,\ \pd M < \infty \}$,

\item $\Fin(A) = \sup\{ \pd(M) |\ M \in \Mod A,\ \pd M < \infty \}$.
\end{itemize}

In the sixties, H. Bass proposes the following two questions about these dimensions, The Finitistic Dimension Conjectures. 

\begin{itemize}

\item[I.] $\fin (A) = \Fin (A)$,

\item[II.] $\fin(A)< \infty$.

\end{itemize}

It is well known that the first one fails, even for monomial algebras (see \cite{ZH}). However the second one,  nowadays called the (Small) Finitistic Dimension Conjecture (FDC), is still an open question. 


Igusa-Todorov algebras are introduced by Wei in \cite{W1} based on the work of Xi (see \cite{X1} and \cite{X2}). In this article, Wei proves that the Finitistic Dimension Conjecture holds for Igusa-Todorov algebras using the Igusa-Todorov functions (see \cite{IT}). 
However, T. Conde shows in her PhD thesis (see \cite{C}) that the exterior algebras $\Lambda(\Bbbk^m)$ are not Igusa-Todorov algebras for $\Bbbk$ an uncontable field and $m\geq 3$. This is a consequence of a Rouquier's result from \cite{R}. 

Lat-Igusa-Todorov algebras (LIT algebras), introduced by Bravo, Lanzilotta, Mendoza and Vivero in \cite{BLMV}, are a generalization of Igusa-Todorov algebras and they also verify FDC. However, not every Artin algebra is a LIT algebra (see \cite{BM}), then the following is a natural question: If $A$ is a LIT algebra, then $A^{op}$ is also LIT? This question was stated by Vivero in \cite{V}.\\


Let $A$ and $B$ be two Artin algebras, $Y$ an $A$-$B$-bimodule, $X$ a $B$-$A$-bimodule,
$\alpha : X \otimes_A Y \rightarrow B$ a $B$-$B$-bimodule homomorphism, and $\beta : Y \otimes_B X \rightarrow A$ an $A$-$A$-bimodule homomorphism. Then from the Morita context $M = (A, Y, X, B, \alpha, \beta)$ we define the Morita context algebra:
$$ \Lambda_{(\alpha, \beta)} = \begin{pmatrix}
A & Y\\
X & B
\end{pmatrix},$$
where the multiplication in $\Lambda_{(\alpha, \beta)}$ is given by 
$$\begin{pmatrix}
a & y\\
x & b
\end{pmatrix} \cdot \begin{pmatrix}
a' & y'\\
x' & b'
\end{pmatrix} = \begin{pmatrix}
aa'+\beta(y\otimes x') & ay'+yb'\\
xa' + bx' & bb' + \alpha(x \otimes y')
\end{pmatrix}$$
and the maps $\alpha$ and $\beta$ satisfy $\alpha(x \otimes y)x' = x\beta(y \otimes x' )$ and $y\alpha(x \otimes y') = \beta(y \otimes x)y' $, to make $\Lambda_{(\alpha, \beta)}$ associative.

We say that a Morita context algebra has zero bimodule homomorphisms if $\alpha = \beta = 0$.

The family of Morita context algebras with zero bimodule homomorphisms has been extensively studied from the homological point of view. We can cite for example, \cite{BM2}, \cite{CRS}, \cite{GP*} and \cite{GP}.  \\


In this article, we study Igusa-Todorov and LIT algebras on Mortita context algebras under similar conditions as \cite{BM2}. The hypotheses are the following:
If $A = \frac{\Bbbk Q_A}{I_A}$ and $B = \frac{\Bbbk Q_B}{I_B}$ are finite dimensional algebras, $C = \frac{\Bbbk \Gamma}{I_C}$ verifies

\begin{itemize}

\item[{\bf H1}:] ${Q_C}_0 = {Q_A}_0 \cup {Q_B}_0$.

\item[{\bf H2}:] ${Q_C}_1 = {Q_A}_1 \cup {Q_B}_1 \cup \{\alpha_j : \start(\alpha_j) \in {Q_A}_0, \ \target(\alpha_j) \in {Q_B}_0 \}_{j \in J} \cup  \{\beta_k : \start(\beta_k) \in {Q_B}_0, \ \target(\beta_k) \in {Q_A}_0 \}_{k \in K}$.

\item[{\bf H3}:] $\langle I_A, I_B, \alpha \alpha_j, \beta \beta_k \text{ for }\alpha \in {Q_A}_1, \beta\in {Q_B}_1, \alpha_j\beta_k, \beta_k\alpha_j \text{ where } j\in J, k\in K\rangle \subset I_C$.

\item[{\bf H4}:] $\mathcal{O}= \langle \add \Orb_{\Omega_A}(\Pi_A(\Omega_C(B_0))) \times \add \Orb_{\Omega_B}( \Pi_B(\Omega_C(A_0))) \rangle \subset K_0(C)$ is finitely generated.

\end{itemize}

We recall from \cite{BM2} that: 

\begin{itemize}

\item Hypotheses {\bf H1}, {\bf H2} and {\bf H3} imply that for all $M \in \mod C$, $\Omega(M) = N_A \oplus N_B$ where $Y_A \in \mod A$ and $Y_B \in \mod B$. They also imply that the bimodules are left-semisimple.

\item Hypotheses {\bf H1}, {\bf H2} and {\bf H3} with an equality in Hypothesis {\bf H3} imply that the bimodules $X$ and $Y$ are right-projective.

\item {\bf H4} is straightforward if $A$ and $B$ are syzygy finite or in case of the inclusion in 
Hypothesis {\bf H3} is an equality.

\end{itemize}

The content of this article can be summarised as follows.
Section 2 is devoted to collecting some necessary material for the development of this work.
We prove, in Section 3, that Morita context algebras which verify {\bf H1}, {\bf H2}, {\bf H3} and {\bf H4}, and arise from Igusa-Todorov algebras (syzygy finite algebras) are Igusa-Todorov algebras (syzygy finite algebras). We prove that the opposite algebra is also Igusa-Todorov (LIT) if the Morita context algebra verifies {\bf H1}, {\bf H2} and {\bf H3} with equallity in {\bf H3}. In section 4, we prove that $\phi_A^{-1}(0)$ is a 0-Igusa-Todorov subcategory if and only if $A$ is a selfinjective algebra or it has finite global dimension. 
Finally, in Section 5, we exhibit an example of an Artin algebra $A$ that is a LIT algebra, but its opposite algebra $A^{op}$ is not.

\section{Preliminaries}

Throughout this article, $A$ is an Artin algebra and $\mod A$ is the category of finitely generated right $A$-modules, $\ind A$ is the subcategory of $\mod A$ formed by all indecomposable modules, $\mathcal{P}_A \subset \mod A$ is the class of projective $A$-modules and $\mathcal{GP}_A \subset \mod A$ is the class of Gorenstein projective $A$-modules. $\mathcal{S} (A)$ is the set of isoclasses of simple $A$-modules and $A_0 = \oplus_{S \in \mathcal{S}(A)}S$. For $M\in \mod A$ we denote by $M^k = \oplus_{i=1}^k M$, by $P(M)$ its projective cover, by $\Omega(M)$ its syzygy and by $\pd (M)$ its projective dimension. The class $\mathcal{P}^{< \infty}$ is formed by all the modules with finite projective dimension. 
The set $\Orb_{\Omega(M)}$ is the $\Omega$-orbit of the module $M$, i.e. $\Orb_{\Omega(M)} = \{\Omega^n(M)\}_{n\geq 0}$. For a subcategory $\mathcal{C} \subset \mod A$, we denote by $\findim (\mathcal{C})$, $\gl (\mathcal{C})$ its finitistic dimension and its global dimension respectively and by $\add \mathcal{C}$ the full subcategory of $\mod A$ formed by all the sums of direct summands of every $M \in \mathcal{C}$.

\subsection{Quivers and path algebras}

If $Q = (Q_0,Q_1,\start,\target)$ is a finite connected quiver, $\Bbbk Q$ denotes its associated {\bf path algebra}. We compose paths in $Q$ from left to right. 
Given $\rho$ a path in $\Bbbk Q$, $\length(\rho)$, $\start(\rho)$ and $\target(\rho )$ denote the length, start and target of $\rho$ respectively. For a quiver $Q$, we denote by $J_Q$ the ideal of $\Bbbk Q$ generated by all its arrows. If $I$ is an admissible ideal of $\Bbbk Q$ ($J^2_Q \subset I \subset J^m_Q$ for some $m \geq 2$), we say that $(Q,I)$ is a {\bf bounded quiver} and the quotient algebra $\frac{\Bbbk Q}{I}$ is the {\bf bound quiver algebra}. A {\bf relation} $\rho$ is an element in $\Bbbk Q$ such that $\rho = \sum \lambda_i w_i$ where the $\lambda_i$ are scalars (not all zero) and the $w_i$ are paths with $\length(w_i) \geq 2$ such that $\start(w_i) = \start(w_j)$ and $\target(w_i) = \target (w_j)$ if $i \not = j$. 
We recall that an admissible ideal $I$ is always generated by a finite set of relations (for a proof see Chapter II.2 Corollary 2.9 of \cite{ASS}).

For a quiver $Q$, we say that $M = (M^v, T_{\alpha})_{v\in Q_0, \alpha \in Q_1}$ is a {\bf representation} if
\begin{itemize}
\item $M^v$ is a $\Bbbk$-vectorial space for every $v\in Q_0$,
\item $T_{\alpha} :M^{\start(\alpha)} \rightarrow M^{\target(\alpha)}$ for every $\alpha \in Q_1$.
\end{itemize}

A representation $(M^v, T_{\alpha})_{v\in Q_0, \alpha \in Q_1}$ is finite dimensional if $M^v$ is finite dimensional for every $v \in Q_0$.

For a path $w = \alpha_1\ldots \alpha_n$ we define $T_{w} = T_{\alpha_1} \ldots T_{\alpha_n}$, and for a relation $\rho = \sum \lambda_i w_i $ we define $T_{\rho} = \sum \lambda_i T_{w_i}$. 

A representation $M = (M^v, T_{\alpha})_{v\in Q_0, \alpha \in Q_1}$ of $Q$ is {\bf bound by} $I$ if we have  $T_{\rho} = 0$ for all relations $\rho \in I$.

Let $M= (M^v, T_{\alpha})_{v\in Q_0, \alpha \in Q_1}$ and $\bar{M} = (\bar{M}^v, \bar{T}_{\alpha})_{v\in Q_0, \alpha \in Q_1}$ be two representations of the bounded quiver $(Q, I)$, a {\bf morphism} $f : M \rightarrow \bar{M}$ is a family $f=(f_v)_{v\in Q_0}$ of $\Bbbk$-linear maps $f_v: M^v \rightarrow \bar{M}^v$ such that for all arrow $\alpha : v \rightarrow w$ we have the following commutative diagram.

$$\xymatrix{ M^v \ar[r]^{T_{\alpha}} \ar[d]_{f_v} & M^w \ar[d]^{f_w}  \\ \bar{M}^v \ar[r]^{\bar{T}_{\alpha}}& \bar{M}^w}$$ 

We denote by $\Rep_{\Bbbk} (Q,I)$ the category of representations of $(Q, I)$ and by $\rep_{\Bbbk} (Q,I)$ the subcategory of $\Rep_{\Bbbk} (Q,I)$ formed by finite dimensional representations of $(Q, I)$. 

\begin{defi}
For an algebra $A = \frac{\Bbbk Q}{I}$ we define the quiver $Q^{\infty}$ as the full subquiver $Q^{\infty} \subset Q$ where a vertex $v \in Q^{\infty}(A)$ if the simple module $S_v$ has $\pd S_v = \infty$.
\end{defi}

We recall that there is a $\Bbbk$-linear equivalence of categories $$F: \Rep_{\Bbbk} (Q,I) \rightarrow \Mod \frac{\Bbbk Q}{I}$$ that restricts to an equivalence of categories $G: \rep_{\Bbbk} (Q,I) \rightarrow \mod \frac{\Bbbk Q}{I}$ (for a proof see Chapter III.1, Theorem 1.6 of \cite{ASS}).

If $A = \frac{\Bbbk Q_A}{I_A}$ and $B = \frac{\Bbbk Q_B}{I_B}$ are finite dimensional algebras, and $C = \frac{\Bbbk \Gamma}{I_C}$ the Morita context algebra has the conditions {\bf H1} and {\bf H2}, then the functors
$\prod_A: \mod C \rightarrow \mod A$ and $\prod_B: \mod C \rightarrow \mod B$ are restictions of the representations, i.e. For a representation $M=(M^v, T_{\alpha})_{v\in {Q_C}_0, \alpha \in {Q_C}_1} \in \Rep_{\Bbbk}(Q_C,I_C)$, we have 

\begin{itemize}
\item $\prod_A(M) = (M^v, T_{\alpha})_{v\in {Q_A}_0, \alpha \in {Q_A}_1} $ and
\item  $\prod_B(M) = (M^v, T_{\alpha})_{v\in {Q_B}_0, \alpha \in {Q_B}_1}$. 
\end{itemize}
For a morphism $f = (f_v)_{v\in {Q_C}_0}: M\rightarrow \bar{M}$, we have 

\begin{itemize}
\item $\prod_A(f) = (f_v)_{v\in {Q_A}_0}: \prod_A(M) \rightarrow \prod_A(\bar{M})$ and
\item  $\prod_B(f) = (f_v)_{v\in {Q_B}_0}: \prod_B(M) \rightarrow \prod_B(\bar{M})$.
\end{itemize}

We add the following lemma from \cite{BM2} since it helps to compute the syzygies in Morita context algebras under hypotheses {\bf H1}, {\bf H2}, {\bf H3} and {\bf H4}.  

\begin{lema}\label{lemita}
Let $A = \frac{\Bbbk Q_A}{I_A}$ and $B = \frac{\Bbbk Q_B}{I_B}$ be finite dimensional algebras. Consider $C = \frac{\Bbbk Q_C}{I_C}$ with the following conditions:

\begin{itemize}

\item ${Q_C}_0 = {Q_A}_0 \cup {Q_B}_0$.

\item ${Q_C}_1 = {Q_A}_1 \cup {Q_B}_1 \cup \{\alpha_j\ \vert \ \start(\alpha_j) \in {Q_A}_0, \ \target(\alpha_j) \in {Q_B}_0 \}_{j \in J} \cup  \{\beta_j\ \vert \ \start(\beta_k) \in {Q_B}_0,\ \target(\beta_k) \in {Q_A}_0 \}_{k \in K}$.

\item $\langle I_A, I_B, \alpha \alpha_j, \beta \beta_k \text{ for }\alpha \in {Q_A}_1, \beta\in {Q_B}_1, \alpha_j\beta_k, \beta_k\alpha_j \text{ where } j\in J, k\in K\rangle \subset I_C$.\\

Then for every module $M \in \mod C$, $\Omega_C(M) = M_1 \oplus M_2$ with $M_1 \in \mod A$ and $M_2 \in \mod B$. Suppose $\Omega_C(\topp(M))= \bar{M}_1 \oplus \bar{M}_2$ with $\bar{M}_1 \in \mod A$ and $\bar{M}_2 \in \mod B$, in case $M \in \mod A$, then $\bar{M}_2 = M_2$, and in case $M \in \mod B$, then $\bar{M}_1 = M_1$. 
\end{itemize}

\end{lema}

\subsection{Igusa-Todorov functions}

In this section, we exhibit some general facts about the Igusa-Todorov functions for an Artin algebra A.
The aim is to introduce material which will be used in the following sections.

\begin{lema}(Fitting Lemma)
Let $R$ be a noetherian ring. Consider a left $R$-module $M$ and $f \in \enn_R (M)$. Then, for any finitely generated $R$-submodule $X$ of $M$, there is a non-negative integer
$$\eta_{f}(X)= \min\{ k \text{ a non-negative integer}: f \vert _{f^m (X)} : f^m (X) \rightarrow f^{m+1} (X), \text{ is injective } \forall m \geq k\}.$$
Furthermore, for any $R$-submodule $Y$ of $X$, we have that $\eta_f(Y) \leq \eta_f(X)$.
\end{lema}

\begin{defi}(\cite{IT})
Let  $K_0(A)$ be the abelian group generated by all symbols $[M]$, with $M \in \mod A$, modulo the relations
\begin{enumerate}
  \item $[M]-[M']-[M'']$ if  $M \cong M' \oplus M''$,
  \item $[P]$ for each projective module $P$.
\end{enumerate}
\end{defi}

Let $\bar{\Omega}: K_0 (A) \rightarrow K_0 (A)$ be the group endomorphism induced by $\Omega$, i.e. $\bar{\Omega}([M]) = [\Omega(M)]$. We denote by $K_i (A) = \bar{\Omega}(K_{i-1}(A))= \ldots = \bar{\Omega}^{i}(K_{0} (A))$ for $i \geq 1$. 
We say that a subgroup $G \subset K_0(A)$ is syzygy finite if there is $n \geq 0$ such that $\bar{\Omega}^n(G)$ is finitely generated. 
For $M\in \mod A$, $\langle \add M\rangle$ denotes the subgroup of $K_0 (A)$ generated by the classes of indecomposable summands of $M$.

 For a subcategory $\mathcal{C} \subset \mod A$, we denote by $\langle\mathcal{C}\rangle \subset K_0(A)$ the free abelian group generated by the classes of direct summands of modules of $\mathcal{C}$.

\begin{defi}\label{monomorfismo}(\cite{IT}) 
The \textbf{(right) Igusa-Todorov function} $\phi$ of $M\in \mod A$  is defined as 
\[\phi_{A}(M) = \eta_{\bar{\Omega}}(\langle \add M \rangle)\]
 In case there is no possible misinterpretation, we will use the notation $\phi$ for the Igusa-Todorov $\phi$ function.
\end{defi}

\begin{prop}(\cite{HLM} \cite{IT}) \label{it1}
If $M,N\in\mod A$, then we have the following.

\begin{enumerate}
  \item $\phi(M) = \pd (M)$ if $\pd (M) < \infty$.
  \item $\phi(M) = 0$ if $M \in \ind A$ and $\pd(M) = \infty$.
  \item $\phi(M) \leq \phi(M \oplus N)$.
  \item $\phi\left(M^{k}\right) = \phi(M)$ for $k \geq 1$.
  \item $\phi(M) \leq \phi(\Omega(M))+1$.
\end{enumerate}

\end{prop}

\begin{prop}\label{invariante}
Suppose $G \subset K_0(A)$ is a finitely generated subgroup with $\rk(G) = m$.
\begin{enumerate}
\item If $\bar{\Omega}(G) \subset G$, then $\eta_{\bar{\Omega}}{}_{\vert G} \leq m$.

\item If $G$ is syzygy finite, then $\eta_{\bar{\Omega}}{}_{\vert G} < \infty$

\end{enumerate}
\begin{proof}
The proof of item 1 is similar to the proof of Proposition 3.6 (item 3) from \cite{LMM}.
The proof of item 2 is similar to the proof of Theorem 3.2. from  \cite{LMM}
\end{proof}
\end{prop}

The result below follows directly from the fact that the Igusa-Todorov function verifies 
$$\phi(M) = \min\{l:\Omega{\vert}_{\Omega^{l+s} \langle \add M \rangle} \text{ is a monomorphism } \forall s \in \mathbb{N}\}$$

\begin{prop}
Given $M \in  \mod A$, 
$$\phi(A) = \max\{n\in \mathbb{N}: \bar{\Omega}^n(v)=0\text{ and }\bar{\Omega}^{n-1}(v)\not=0 \text{ for some } v\in \langle \add M \rangle\}.$$

\end{prop}

\subsection{LIT algebras}

Lat-Igusa-Todorov algebras were introduced in \cite{BLMV} as a generalization of Igusa-Todorov algebras (see Definition 2.2 of \cite{W1}). They also verify the Finitistic Dimension Conjecture as can be seen in Theorem 5.4 of \cite{BLMV}. 

Let $A$ be an Artin algebra. If $\mathcal{D} \subset \mod A$ is a subcategory such that
\begin{enumerate}
\item $\mathcal{D} = \add (\mathcal{D})$,

\item  $\Omega (\mathcal{D}) \subset \mathcal{D}$ and

\item  $\fidim(\mathcal{D}) = 0$,
\end{enumerate}
we call it a {\bf $0$-Igusa-Todorov subcategory}.

\begin{ej}\label{clases LIT}
Let $A$ be an Artin algebra.
\begin{enumerate}
\item If $\fidim (A) = 0$, then $\mathcal{D} = \mod A$ is a $0$-Igusa-Todorov subcategory.

\item If $\fidim (A) = 1$, then $\mathcal{D} = \Omega (\mod A)$ is a $0$-Igusa-Todorov subcategory.

\item Note that the two previous items cannot be generalized to algebras with $\fidim (A) = n > 1$.  On one hand $\Omega^{n}(\mod A)$ is not closed by direct summands in general and on the other hand $\phi(\add \Omega^{n}(\mod A))$ is not $0$ generally. 

\item $\mathcal{G}\mathcal{P}(A)$ and $^\bot A$ are $0$-Igusa-Todorov subcategories.
\end{enumerate}
\end{ej}

\begin{defi}(\cite{BLMV})\label{modulo LIT}
Let $A$ be an Artin algebra. A subcategory $ \mathcal{C} \subset \mod A$ is called ${\bf (n, V, \mathcal{D}){\text -}}${\bf Lat-Igusa-Todorov} (for short ${\bf (n, V, \mathcal{D})\text{-}LIT}$) if the following conditions are verified 
\begin{itemize}

\item There is some $0$-Igusa-Todorov subcategory $\mathcal{D} \subset \mod A$,

\item there is some $V \in \mod A$ satisfying that each $M \in \mathcal{C}$ admits an exact sequence:
$$\xymatrix{0 \ar[r]& V_1 \oplus D_1 \ar[r] & V_0 \oplus D_0 \ar[r] & \Omega^n(M)\ar[r] & 0}$$

such that $V_0, V_1 \in \add (V)$ and $D_0, D_1 \in \mathcal{D}$.
\end{itemize}
We say that $V$ is an ${\bf (n, V, \mathcal{D})\text{-}}${\bf Lat-Igusa-Todorov} {\bf module} (for short a ${\bf n\text{-}LIT}$ {\bf module}) for $\mathcal{C}$.
\end{defi}

\begin{defi}(\cite{BLMV})\label{algebra LIT}
We say that $A$ is an  ${\bf (n, V, \mathcal{D})\text{-}}${\bf Lat-Igusa-Todorov} {\bf algebra} (for short an  ${\bf  (n, V, \mathcal{D})\text{-}LIT}$ {\bf algebra}) if $\mod A$ is $(n, V, \mathcal{D})\text{-}${\rm{LIT}}. We say that $A$ is a {\rm{LIT}} algebra if $A$ is $n\text{-}${\rm{LIT}} for some non-negative integer $n$.
\end{defi}

LIT algebras are a generalization of Igusa-Todorov algebras which are defined by Wei in \cite{W1}. 

\begin{defi} (\cite{W1})
If $\mathcal{D} = \{0\}$ in Definition \ref{algebra LIT}, we say that $A$ is an {\bf n-Igusa-Todorov algebra} with $V$ its {\bf n-Igusa-Todorov module}.
\end{defi}

The following result shows the interest in studying LIT algebras. These algebras verify the finitistic dimension conjecture. For proof see \cite{BLMV}.

\begin{teo}(\cite{BLMV})\label{LIT finitista} Let $A$ be an $(n, V, \mathcal{D})$-{\rm{LIT}} algebra. Then
$$\fin (A) \leq \psi_{[\mathcal{D}]} (V) + n + 1 < \infty.$$
\end{teo}

The following two results build new families of Igusa-Todorov and LIT algebras over Morita context algebras respectively.  
The first one is proved by Yajun Ma and Nanqing Ding in \cite{DM}. The second result is proved by Vivero in \cite{V}, and we show in Example \ref{ejemplin} that item 1 of this theorem is not true for a less restrictive hypothesis.


\begin{teo}(\cite{DM})\label{DM}
Let $A_{(0,0)} =  \left( \begin{array}{cc} T & N \\ M & U  \end{array} \right)$ be a Morita ring such that $_A N, N_B M_A , _BM$ are projective modules. Then for any positive integer n, we have

\begin{enumerate}

\item $A_{(0,0)}$ is n-syzygy-finite if and only if T and U are n-syzygy-finite.

\item $A_{(0,0)}$ is n-Igusa–Todorov if and only if T and U are n-Igusa-Todorov.

\end{enumerate}

\end{teo}

\begin{teo}(\cite{V})\label{Algebra LIT triangular} Let $T$ and $U$ be $(n, V_T , D_T )$ and $(n, V_U , D_U )$ LIT algebras. Let $ _UM_T$
be projective both as a left $U$-module and as a right $T$-module and such that
$M \otimes_T$ $P$ is indecomposable whenever $P \in \mathcal{P}_T$ is indecomposable. Then

\begin{enumerate}

\item The algebra $A = \left( \begin{array}{cc} T & 0 \\ M & U  \end{array} \right)$
is $(n + 1, V_A , D_A )$-LIT, where the class $D_A = \add (\Omega (D_T , 0, 0) \oplus (0, D_U , 0))$
and $V_A = \Omega(V_T , V_U , 0)  \oplus A$.

\item $\findim(A) < \infty$.

\end{enumerate}

\end{teo} 


\section{Morita context algebras and Igusa-Todorov algebras}

In this sections, we give similar results to Theorems \ref{DM} and \ref{Algebra LIT triangular} (Propositions \ref{sizigia finita cuadrada}, \ref{Igusa-Todorov cuadrada} and \ref{C LIT si A es IT y B LIT}), changing the hypothesis on the modules $M$ and $N$. We also give a dual result (Proposition \ref{dual}) with a more restrictive hypothesis.

Consider $A$, $B$ and $C$ that verify {\bf H1}, {\bf H2} and {\bf H3}. Recall that $\mathcal{O}= \langle \add \Orb_{\Omega_A}(\Pi_A(\Omega_C(B_0))) \times \add \Orb_{\Omega_B}( \Pi_B(\Omega_C(A_0))) \rangle \subset K_0(C)$. Let $\mathcal{A}$ and $\mathcal{B}$ be the subgroups of $K_1(C)$ generated by the classes of indecomposable $A$-modules and $B$-modules respectively which do not belong to $\mathcal{O}$. 

The following result is a generalization of Item 2 of Remark 3.6 from \cite{BM2}.

\begin{prop}\label{sizigia finita cuadrada}

Let $A = \frac{\Bbbk Q_A}{I_A}$ and $B = \frac{\Bbbk Q_B}{I_B}$ be finite dimensional algebras. Consider $C = \frac{\Bbbk \Gamma}{I_C}$ with the following conditions:

\begin{itemize}

\item ${Q_C}_0 = {Q_A}_0 \cup {Q_B}_0$.

\item ${Q_C}_1 = {Q_A}_1 \cup {Q_B}_1 \cup \{\alpha_j : {Q_A}_0 \rightarrow {Q_B}_0 \}_{j \in J} \cup  \{\beta_k : {Q_B}_0 \rightarrow {Q_A}_0 \}_{k \in K}$.

\item $\langle I_A, I_B, \alpha \alpha_j, \beta \beta_k \text{ for }\alpha \in {Q_A}_1, \beta\in {Q_B}_1, \alpha_j\beta_k, \beta_k\alpha_j \text{ where } j\in J, k\in K\rangle \subset I_C$.

\item $\mathcal{O}=\langle \add \Orb_{\Omega_A}(\Pi_A(\Omega_C(C_0))) \times \add \Orb_{\Omega_B}( \Pi_B(\Omega_C(C_0)))\rangle \subset K_0(C)$ is finitely generated.

\end{itemize}

Then $C$ is a syzygy finite algebra if and only if $A$ and $B$ are syzygy finite algebras. 
\begin{proof}

Suppose $A$ and $B$ are finite syzygy, in particular we can assume $K_n(A) \subset \langle [N_1],[N_2],\ldots,[N_t]\rangle$ and $K_n(B) \subset \langle [\tilde{N}_1],[\tilde{N}_2],\ldots,[\tilde{N}_s]\rangle$. We also assume that $\mathcal{O} = \langle [O_1], [O_2], \ldots,[O_m] \rangle$.

By Lemma \ref{lemita}, we have that $K_1(C) \subset \mathcal{A} \times \mathcal{B} \times  \mathcal{O}$, and $$\bar{\Omega}_C ([M]) = \left\lbrace \begin{array}{ll} \bar{\Omega}_A([M]) + [O] & \text{ if } [M]\in \mathcal{A}, \\  \bar{\Omega}_B([M]) + [O]  & \text{ if } [M]\in \mathcal{B}, \\ {} [N] & \text{ if } [M]\in \mathcal{O}, \end{array} \right.$$ 
where $O$ belongs to $\mathcal{O}$.
 
Hence if $[M] \in K_1(C)$, we can consider three possible cases.

\begin{itemize}
\item $[M] \in \mathcal{O}$,

\item $[M] \in \mathcal{B}$,

\item $[M] \in \mathcal{A}$.

\end{itemize}

In the first case we have that $\bar{\Omega}_C^{k'}([M]) \subset \langle [O_1], [O_2], \ldots,[O_m] \rangle$ $\forall k' \geq k$. In the second case, we know that $\bar{\Omega}_C^n([M]) = \sum_{i=1}^t\alpha_i[\tilde{N}_i] + [D]$ with $[D] \in \mathcal{O}$, then $\bar{\Omega}_C^{n+k}([M]) = \sum_{i=1}^s\alpha_i[\bar{\Omega}_C^{k}(\tilde{N}_i)] + \sum_{i=1}^m\beta_i[O_i]$. Analogously $\bar{\Omega}_C^{n+k}([M]) = \sum_{i=1}^t\alpha_i[\bar{\Omega}_C^{k}({N}_i)] + \sum_{i=1}^m\beta_i[O_i]$ for $[M] \in \mathcal{A}$. We finally deduce that 
$$K_{k+n}(C) \subset \langle \add ((\oplus_{i=1}^m O_i )\oplus (\oplus_{i=1}^t\bar{\Omega}_C^{k}({N}_i)) \oplus (\oplus_{i=1}^s\bar{\Omega}_C^{k}(\tilde{N}_i))) \rangle.$$
\end{proof}

\end{prop}

\begin{prop}\label{Igusa-Todorov cuadrada}

Let $A = \frac{\Bbbk Q_A}{I_A}$ and $B = \frac{\Bbbk Q_B}{I_B}$ be finite dimensional algebras. Consider $C = \frac{\Bbbk \Gamma}{I_C}$ with the following conditions:

\begin{itemize}

\item ${Q_C}_0 = {Q_A}_0 \cup {Q_B}_0$.

\item ${Q_C}_1 = {Q_A}_1 \cup {Q_B}_1 \cup \{\alpha_j : {Q_A}_0 \rightarrow {Q_B}_0 \}_{j \in J} \cup  \{\beta_k : {Q_B}_0 \rightarrow {Q_A}_0 \}_{k \in K}$.

\item $\langle I_A, I_B, \alpha \alpha_j, \beta \beta_k \text{ for }\alpha \in {Q_A}_1, \beta\in {Q_B}_1, \alpha_j\beta_k, \beta_k\alpha_j \text{ where } j\in J, k\in K\rangle \subset I_C$.

\item $\mathcal{O}= \langle\add \Orb_{\Omega_A}(\Pi_A(\Omega_C(C_0))) \times \add \Orb_{\Omega_B}( \Pi_B(\Omega_C(C_0)))\rangle \subset K_0(C)$ is finitely generated.

\end{itemize}

Then $C$ is an (n+1)-Igusa-Todorov algebra if $A$ and $B$ are n-Igusa-Todorov algebras. 
\begin{proof}

Suppose that $n$ is the minimum natural number such that $A$ and $B$ are n-Igusa-Todorov simultaneously.

We know that $K_1(C) \subset \mathcal{A} \times \mathcal{B} \times  \mathcal{O}$. Since $A$ and $B$ are $n$-Igusa-Todorov algebras there exist a $A$-module $V_A$ and a $B$-module $W_B$ such that for all  $ M, N \in \mod C$ where $[M] \in \mathcal{A}$ and $[N] \in \mathcal{B}$, there are two short exact sequences
$$\xymatrix{0\ar[r] & V_A'\ar[r] & V_A''\ar[r] & \Omega^n_A(M) \ar[r] & 0} \text{ and } \xymatrix{0\ar[r] & W_B'\ar[r] & W_B''\ar[r] & \Omega^n_B(N) \ar[r] & 0},$$
with $V_A', V_A'' \in \add V_A$ and  $W_B', W_B'' \in \add W_B$.

Consider $M \in \Omega_C^{n+1}(\mod C)$, then $[M] = [M_0]+[M_1]+[M_2]$ with $[M_0] \in \mathcal{O}$, $[M_1] \in \bar{\Omega}_A^{n}(\mathcal{A})$ and $[M_2] \in \bar{\Omega}_B^{n}(\mathcal{B})$.

We conclude that $C$ is a $(n+1)$-Igusa-Todorov algebra and $V_A\oplus A \oplus W_B \oplus B \oplus (\oplus_{O\in I} O)$
a Igusa-Todorov module with $I= \{ O \in \ind C: [O] \in \mathcal{O}\}$. 
\end{proof}
\end{prop}

\begin{prop}\label{C LIT si A es IT y B LIT}
Let $A = \frac{\Bbbk Q_A}{I_A}$ and $B = \frac{\Bbbk Q_B}{I_B}$ be finite dimensional algebras. Consider $C = \frac{\Bbbk \Gamma}{I_C}$ with the following conditions:

\begin{itemize}

\item ${Q_C}_0 = {Q_A}_0 \cup {Q_B}_0$.

\item ${Q_C}_1 = {Q_A}_1 \cup {Q_B}_1 \cup \{\alpha_j : {Q_A}_0 \rightarrow {Q_B}_0 \}_{j \in J}$.

\item $\langle I_A, I_B, \alpha \alpha_j \text{ for }\alpha \in {Q_A}_1,  \text{ where } j\in J \rangle \subset I_C$.

\item $\mathcal{O} = \langle \add \Orb_{\Omega_B}( \Pi_B(\Omega_C(C_0)))\rangle \subset K_0(C)$ is finitely generated.

\end{itemize}

Then $C$ is an (n+1)-LIT algebra if $A$ is an n-Igusa-Todorov algebra and $B$ is an n-LIT algebra. 

\begin{proof}
Suppose that $A$ is an $n$-Igusa-Todorov algebra with an Igusa-Todorov module $W$ and $B$ is an $(n, V, \mathcal{D})$-LIT algebra.

Observe that for $M\in \mod C$, $\Omega_C(M) = M_1 \oplus M_2$ where $M_1 \in \mod A$ and $M_2 \in \mod B$, and
$$\Omega_C(N) = \left\lbrace \begin{array}{ll} \Omega_A(N) \oplus O & \text{ if } N\in \mod A, \\ \Omega_B(N)& \text{ if }  N\in \mod B,\end{array} \right.$$
where $O$ belongs to $\mathcal{O}$.

Then $C$ is an $(n+1, V\oplus W \oplus A \oplus (\oplus_{O \in I} O), \mathcal{D})$-LIT algebra and $V\oplus W \oplus A \oplus (\oplus_{O \in I} O)$ is a $n$-LIT module where $I = \{O \in \ind C: [O] \in \mathcal{O}\}$
\end{proof}

\end{prop}

With a more restrictive hypothesis on the ideal $I_C$ we can consider a weaker hypothesis on the previous result.

\begin{prop}
Let $A = \frac{\Bbbk Q_A}{I_A}$ and $B = \frac{\Bbbk Q_B}{I_B}$ be finite dimensional algebras. Consider $C = \frac{\Bbbk \Gamma}{I_C}$ with the following conditions:

\begin{itemize}

\item ${Q_C}_0 = {Q_A}_0 \cup {Q_B}_0$.

\item ${Q_C}_1 = {Q_A}_1 \cup {Q_B}_1 \cup \{\alpha_j : {Q_A}_0 \rightarrow {Q_B}_0 \}_{j \in J} \cup  \{\beta_k : {Q_B}_0 \rightarrow {Q_A}_0 \}_{k \in K}$, such that $\start(\alpha_j) \not = \target(\beta_k)$ and $ \start(\beta_k) \not = \target(\alpha_j)$ $\forall j \in J, \forall k \in K$.

\item $\langle I_A, I_B, \alpha \alpha_j, \beta \beta_k \text{ for }\alpha \in {Q_A}_1, \beta\in {Q_B}_1, \alpha_j\beta_k, \beta_k\alpha_j \text{ where } j\in J, k\in K\rangle = I_C$.


\end{itemize}

Then $C$ is (n+1)-LIT if $A$ and $B$ are n-LIT.

\begin{proof}
Suppose that $A$ is $(n, V_1, \mathcal{D}_1)$-LIT algebra and $B$ is $(n, V_2, \mathcal{D}_2)$-LIT algebra.
We recall that from the proof of Proposition 3.5 from \cite{BM2} we have $$K_1(C) \subset \langle [M], [N] \ \vert \ M \in \mod A, N \in \mod B \rangle = \langle \mathcal{P}_{\partial A}\rangle \times \langle \mathcal{P}_{\partial B} \rangle \times K_0(A) \times K_0(B).$$

where 
\begin{itemize}
\item $\partial A = \{v\in {Q_A}_0 \ \vert \ \exists \ \alpha \in {Q_C}_1 \text{ with } \start(\alpha) = v \text{ and } \target(\alpha) \in {Q_B}_0\}$ and 
\item $\partial B = \{v\in {Q_B}_0 \ \vert \ \exists \ \beta \in {Q_C}_1 \text{ with } \start(\beta) = v \text{ and } \target(\beta) \in {Q_A}_0\}$.
\end{itemize}
We also have for $M_1 \in \mod A$ $\Omega_C(M_1) = \Omega_A(M_1)\oplus P_B$ and for $M_2 \in \mod B$ $\Omega_C(M_2) = \Omega_B(M_2) \oplus P_A$.  Since  $\start(\alpha_j) \not = \target(\beta_k)$ and $ \start(\beta_k) \not = \target(\alpha_j)$ $\forall j \in J, \forall k \in K$ the modules $P_A$ and $P_B$ are in $\mathcal{P}_C$. 

Consider the full subcategories
\begin{itemize}
\item $\tilde{\mathcal{D}}_1 = \{ M: M \in \mathcal{D}_1 \text{ with no direct summands in } \mathcal{P}_A \}$ and 
\item $\tilde{\mathcal{D}}_2 = \{ M: M \in \mathcal{D}_2 \text{ with no direct summands in } \mathcal{P}_B \}$.
\end{itemize} 
Then it is clear that $C$ is a $(n+1, V, \mathcal{D})$-LIT algebra with $ \mathcal{D} = \{D_1\oplus D_2 \oplus P : D_i \in \tilde{\mathcal{D}}_i \text{ for } i=1,2 \text{ and } P \in \mathcal{P}_C\}$ and $V= V_1\oplus V_2 \oplus (\oplus_{v \in \partial A} {P_v}_A) \oplus (\oplus_{v \in \partial B} {P_v}_B) $ where ${P_v}_A$ $({P_v}_B)$ is the indecomposable projective $A$-module (B-module) associated to the vertex $v$.
\end{proof}

\end{prop}

In the next we will consider $A = \frac{\Bbbk Q_A}{I_A}$, $B = \frac{\Bbbk Q_B}{I_B}$ and $C = \frac{\Bbbk \Gamma}{I_C}$ finite dimensional algebras with the following conditions:

\begin{itemize}

\item ${Q_C}_0 = {Q_A}_0 \cup {Q_B}_0$.

\item ${Q_C}_1 = {Q_A}_1 \cup {Q_B}_1 \cup \{\alpha_j : {Q_A}_0 \rightarrow {Q_B}_0 \}_{j \in J} \cup  \{\beta_k : {Q_B}_0 \rightarrow {Q_A}_0 \}_{k \in K}$.

\item $\langle I_A, I_B, \alpha \alpha_j, \beta \beta_k \text{ for }\alpha \in {Q_A}_1, \beta\in {Q_B}_1, \alpha_j\beta_k, \beta_k\alpha_j \text{ where } j\in J, k\in K\rangle = I_C$.

\end{itemize}

We denote by $\mathcal{T}$ the set $\{ S_{v_0}\in \mathcal{S}(C) \text{ and there is } \beta_k \text{ or } \alpha_j \text{ with } \start(\beta_k) = v_0 \text{ or } \start(\alpha_j) = v_0 \}$.

Following Proposition 3.5 from \cite{BM2}, consider $f: K_1(C^{op}) \rightarrow \langle \mod A^{op} \rangle \times \langle \mod B^{op} \rangle \times \langle C_0 \rangle$ the function defined in the indecomposable modules as follows
$$f([M]) = \left\lbrace \begin{array}{ll} (0,0,[S_{v_0}]) & \text{ if } M = S_{v_0}\in \mathcal{T}, \\ ([\prod_A(M)],0,0) & \text{ if we are not in the first case and } \topp(M) \in \mod A^{op}, \\ (0, [\prod_B(M)],0) & \text{ if we are not in the first case and } \topp(M) \in \mod B^{op}. \end{array} \right.$$

We recall Claim 3) of the proof of Proposition 3.5 from \cite{BM2}. $f: K_1(C^{op}) \rightarrow \langle \mod A^{op} \rangle \times \langle \mod B^{op} \rangle \times \langle C_0 \rangle$ is a monomorphism of groups and $f([\Omega_{C^{op}}(M)]) = ([\Omega_{A^{op}}(M_1)], [\Omega_{B^{op}}(M_2)], 0)$ if $f([M]) = ([M_1], [M_2], 0)$.  

Observe that we complete the previous claim as follows. If $f([M]) = (0,0,[S_A\oplus S_B])$ where $S_A \in  \mod A^{op}$ and  $S_B \in \mod B^{op}$, then $f([\Omega_C(M)]) = ([\Omega_{A^{op}}(S_A)], [\Omega_{B^{op}}(S_B)], [S'])$ where $S' \in C_0$.\\

We also define the functor $G_A : \mod A^{op} \rightarrow \mod C^{op}$ as follows

\begin{itemize}
\item If $M = (M_i, T_{\alpha})_{i \in {Q_{A}}_0, \ \alpha \in  {Q_{A}}_1 }$, then $G_A(M) = \tilde{M} = (\tilde{M}_i, \tilde{T}_{\alpha})_{i \in {Q_{A}}_0, \ \alpha \in  {Q_{A}}_1 }$ where

\begin{itemize}
\item $\tilde{M}_i = \left\{ \begin{array}{ll} M_i & \text{ if } i \in {Q_{A}}_0, \\ \oplus_{h \in {Q_{A}}_0}(\oplus_{ \substack{\\ \gamma \in {Q_{C}}_0 \\  \start(\gamma) = h \\ \target(\gamma) = i }} M_h) & \text{ if } i \in \partial B, \\ 0 & \text{ if } {Q_{B}}_0 \setminus \partial B. \end{array} \right.$

\item $\tilde{T}_{\alpha} = \left\{ \begin{array}{ll} T_{\alpha} & \text{ if } \alpha \in {Q_{A}}_1, \\ ( \oplus_{ \substack{ k \in K \\ k \not = k_0}} 0) \oplus_{k = k_0} Id  &  \text{ if } \alpha = \alpha_{k_0} \text{ with } k_0 \in K, \\ 0 &  \text{ if } \alpha = \beta_j \text{ with } j \in J, \\ 0 & \text{ if } \alpha \in {Q_{B}}_1. \end{array} \right.$

\end{itemize}

\item If $\phi = (\phi_{i})_{i \in {Q_{A}}_0 } : M \rightarrow N$ where $M = (M_i, T_{\alpha})_{i \in {Q_{A}}_0, \ \alpha \in  {Q_{A}}_1 }$ and $N = (N_i, S_{\alpha})_{i \in {Q_{A}}_0, \ \alpha \in  {Q_{A}}_1 }$. Then $G_A(\phi) = \tilde{\phi} = (\tilde{\phi}_{i})_{i \in {Q_{C}}_0 } $ where

$$\tilde{\phi}_i =  \left\{\begin{array}{ll} \phi_i & \text{ if } i \in {Q_{A}}_0, \\ \oplus_{\substack{\alpha_k, k \in K \\ \start(\alpha_k) = j \\ \target(\alpha_k) = i}} \phi_i & \text{ if } i \in \partial B, \\ 0 & \text{ if } i \in {Q_{B}}_0 \setminus \partial B. \end{array} \right.$$

\end{itemize}

We define $G_B : \mod B^{op} \rightarrow \mod C^{op}$ analogously.

\begin{obs} Observe that

\begin{itemize}
\item The functors $G_A$ and $G_B$ are exact, additives and projective preseving.  By Remark 3.1 of \cite{BMR2}, the functors $G_A$ and $G_B$ induce the group morphisms $\bar{G}_A : K_0(A) \rightarrow K_0(C)$ and $\bar{G}_B : K_0(B) \rightarrow K_0(C)$ which commute with the morphisms $\bar{\Omega}$.

\item $\bar{G}_A(f([M)])) = [M]$ if $M \in \mod A^{op}$ and $\bar{G}_B(f([(N)])) = [N]$ if $N \in \mod B^{op}$.  
\end{itemize}
\end{obs}

\begin{prop} \label{dual}
Let $A = \frac{\Bbbk Q_A}{I_A}$ and $B = \frac{\Bbbk Q_B}{I_B}$ be finite dimensional algebras. Consider $C = \frac{\Bbbk \Gamma}{I_C}$ with the following conditions:

\begin{itemize}

\item ${Q_C}_0 = {Q_A}_0 \cup {Q_B}_0$.

\item ${Q_C}_1 = {Q_A}_1 \cup {Q_B}_1 \cup \{\alpha_j : {Q_A}_0 \rightarrow {Q_B}_0 \}_{j \in J} \cup  \{\beta_k : {Q_B}_0 \rightarrow {Q_A}_0 \}_{k \in K}$.

\item $\langle I_A, I_B, \alpha \alpha_j, \beta \beta_k \text{ for }\alpha \in {Q_A}_1, \beta\in {Q_B}_1, \alpha_j\beta_k, \beta_k\alpha_j \text{ where } j\in J, k\in K\rangle = I_C$.


\end{itemize}

Then $C^{op}$ is an (n+1)-Igusa-Todorov ((n+1)-LIT) algebras if $A^{op}$ and $B^{op}$ are n-Igusa-Todorov (n-LIT) algebras.

\begin{proof}

From the previous remark, it is easy to check the thesis in both cases. So we only exhibit the Igusa-Todorov modules and the 0-Igusa-Todorov subcategory. Consider $V_3 = \oplus_{i=1}^n \Omega^{i}(C_0)$.


\begin{itemize}
\item Suppose that $n$ is the minimum integer number such that $A^{op}$ and $B^{op}$ are n-Igusa-Todorov algebras simultaneously with n-Igusa-Todorov modules $V_1$ and $V_2$ respectively.

Consider $V = G_A(V_1) \oplus G_B(V_2) \oplus V_3$, then $C^{op}$ is an (n+1)-Igusa-Todorov algebra with (n+1)-Igusa-Todorov module $V$. 

\item Now suppose that $n$ is the minimum natural number such that $A^{op}$ and $B^{op}$ are $(n, V_1, \mathcal{D}_1)\text{-}LIT$ and $(n, V_2, \mathcal{D}_2)\text{-}LIT$ algebras respectively.

Consider $V = G_A(V_1) \oplus  G_B(V_2) \oplus V_3$, and $\mathcal{D} = \{G_A(M)\oplus G_B(N)\oplus P :M \in \mathcal{D}_1, N \in \mathcal{D}_2, P \in \mathcal{P}(C)\}$. Then $C^{op}$ is an $ (n+1, V, \mathcal{D})\text{-}LIT$ algebra. 

\end{itemize}
\end{proof}

\end{prop}

\section{$\phi_A(0)^{-1}$ as a 0-Igusa-Todorov subcategory}

Igusa-Todorov algebras and LIT algebras are different because of the existence of a 0-Igusa-Todorov subcategory $\mathcal{D}$. 
For this reason, it is interesting to know when a subcategory is 0-Igusa-Todorov.
In section 5 of \cite{BM} the authors show examples of algebras where the 0-Igusa-Todorov subcategories are trivial, i.e. $\mathcal{D} \subset \mathcal{P}(A)$ for every $\mathcal{D}$ a 0-Igusa-Todorov subcategory.  
On the other hand, Example \ref{clases LIT} shows some algebras with nontrivial 0-Igusa-Todorov subcategories.

\begin{defi}
For an Artin algebra $A$, we denote by $\phi_A(0)^{-1}$ the class $\{ M \in \mod A: \phi(M) = 0\}$.
\end{defi}

Since for every 0-Igusa-Todorov subcategory $\mathcal{D}$, we have the inclusion $\mathcal{D} \subset \phi_A(0)^{-1}$, in this section we study when the class $\phi_A(0)^{-1}$ is a 0-Igusa-Todorov subcategory for an Artin algebra $A$. In particular we characterize when $\phi_A(0)^{-1}$ is a 0-Igusa-Todorov subcategory for $A = \frac{\Bbbk Q}{I}$ a bound quiver. 

\begin{obs}\label{observacioncita} For an Artin algebra $A$
\begin{enumerate}
\item $\mathcal{P}_A \subset \phi_A(0)^{-1}$.
\item If $M \in \ind A$ and $\pd(M) = \infty$, then $M \in \phi_A(0)^{-1}$.
\end{enumerate}

\end{obs}

\begin{obs} \label{aditivas}
Let $A$ be an Artin algebra. If $\gl (A) < \infty$ or $A$ is selfinyective, then $\phi_A(0)^{-1}$ is additive.
\begin{proof}
If $\gl (A) < \infty$, then $\phi_A(0)^{-1} = \mathcal{P}(A)$ and $ \mathcal{P}(A)$ is additive. If $A$ is selfinjective, then $\phi_A(0)^{-1} = \mod A$ by the main result from \cite{HL} and $\mod A$ is additive. 
\end{proof}
 
\end{obs}

\begin{lema}\label{1}
Let $A$ be an Artin algebra. If $\phi_A(0)^{-1}$ is additive, then $\Omega(M)$ is indecomposable and $\pd (\Omega(M)) = \infty$ for every indecomposable module $M$ with $\pd(M) = \infty$.
\begin{proof}

Since $\pd(M) = \infty$, then it is well known that $\pd(\Omega(M)) = \infty$. 
Consider $M \in \phi_A(0)^{-1}$ an indecomposable module. Since $\phi_A(0)^{-1}$ is additive, by Item 2 of Remark \ref{observacioncita}, we can suppose that $\Omega(M) = N_1 \oplus N_2$ with $N_1 \in \phi_A(0)^{-1}$ and $N_2 \in  \mathcal{P}^{< \infty} (A)$. Let $M_1$ and $M_2$ be such that $\Omega(M_1) = N_1$ and $\Omega(M_2) = N_2$. We have two different cases.

\begin{enumerate}
\item $N_1$ is indecomposable and $N_2 \not \cong \{0\}$.

\item $N_1$ has at least two indecomposable direct summands. 
\end{enumerate}

\begin{enumerate}

\item Since $N_1$ is an indecomposable module, then $M_1$ can be considered indecomposable and $M_1\in \phi_A(0)^{-1}$. 
Since $N_2 \not \cong \{0\}$, then $M_1 \not \cong M$. It is clear that $\phi(M \oplus M_1)> \pd (N_2)> 0$. On the other hand, since $\phi_A(0)^{-1}$ is additive, $M\oplus M_1 \in \phi_A(0)^{-1}$ and $\phi(M\oplus M_1) = 0$, which is a contradiction. 

\item Suppose $ N_1 = \oplus_{j \in J} {\bar{N}_j}$ is a decomposition into indecomposable modules with $\vert J \vert > 1$. For every $j \in J$ we have $\bar{M}_j \in \phi_A(0)^{-1}$ an indecomposable module such that $\Omega(\bar{M}_j) = {\bar{N}_j}$. Then it is clear that $\phi(M \oplus (\oplus_{j\in J} \bar{M}_j)) > 0$. On the other hand, since $\phi_A(0)^{-1}$ is additive, $M\oplus (\oplus_{j\in J} \bar{M}_j)) \in \phi_A(0)^{-1}$ and $\phi(M\oplus (\oplus_{j\in J} \bar{M}_j)) = 0$, which is also a contradiction.

\end{enumerate}

\end{proof}

\end{lema}

\begin{lema} \label{2}
Let $A$ be an Artin algebra. If $\phi_A(0)^{-1}$ is additive, then the restriction $\bar{\Omega} :  \langle \phi_A(0)^{-1} \rangle \rightarrow  \langle\phi_A(0)^{-1} \rangle$ is injective.

\begin{proof}
Suppose that $\bar{\Omega}$ is not injective. Then there are two non-projective modules $M_1$, $M_2$ in $\phi_A(0)^{-1}$ such that $\bar{\Omega}([M_1]-[M_2]) = 0$. Then there are two short exact sequences
$$0 \rightarrow \tilde{M}  \rightarrow P_1 \rightarrow M_1 \rightarrow 0,$$
$$0 \rightarrow \tilde{M} \rightarrow P_2 \rightarrow M_2 \rightarrow 0.$$ Therefore $\phi(M_1\oplus M_2) \geq 1$. This is a contradiction since $\phi_A(0)^{-1}$ is additive.
   
\end{proof}

\end{lema}

\begin{lema}\label{3}
Let $A$ be an Artin algebra such that $\phi_A(0)^{-1}$ is an additive subcategory and $M \in \phi_A(0)^{-1}$ is an indecomposable non-projective $A$-module. For every $S \in \mathcal{S}_A$ with an epimorphism $M \rightarrow S \rightarrow 0$, $\pd(S) = \infty$. 

\begin{proof}

Suppose $\pd S < \infty$. We have the short exact sequence $$0 \rightarrow K \rightarrow M \rightarrow S \rightarrow 0$$
Since $\phi_A(0)^{-1}$ is additive, then $K \cong K_0 \oplus K_1$ with $K_0 \in \phi_A(0)^{-1}$ and $\pd(K_1) < \infty$. Consider $n = \max\{\pd(S), \pd(K_1)\}$. By the Horseshoe Lemma, there is a projective module P such that $\Omega^n(K_0) \cong \Omega^n(K_0\oplus K_1) \cong \Omega^n(M)\oplus P$. Since $\phi_A(0)^{-1}$ is additive $[K_0] = [M]$, hence $K_0 \cong M$. This implies that $K_1 = 0$ and $S = 0$, which is a contradiction. 
\end{proof}
\end{lema}

\begin{lema}\label{4}
Let $A$ be an Artin algebra. If $\phi_A(0)^{-1}$ is an additive subcategory and every simple module has infinite projective dimensión, then $\fin (A) = 0$.

\begin{proof}

The proof is by induction on the Loewy length. Let $M$ be an indecomposable $A$-module.
\begin{itemize}
\item If $\loewy(M) = 1$, then $M$ is semisimple and $\pd(M) = \infty$.

\item If $\loewy(M) = 2$, then there is a short exact sequence
$$0 \rightarrow S \rightarrow M \rightarrow S' \rightarrow 0 $$
where $S = \soc(M)$ and $S' = \topp(M)$. Then, by the Horseshoe Lemma, $\Omega^{m}(S) = \Omega^{m+1}(S')$. Since $S$ and $S'$ belong to $\phi_A(0)^{-1}$, then $\Omega(S') = S$. Hence we have the following short exact sequence
$$0 \rightarrow S \rightarrow P_{S'} \rightarrow S' \rightarrow 0$$
Since $\topp (M) = S'$, we know that there exists an epimorphism $\pi : P_{S'} \rightarrow M$, and we deduce that $\pi$ is an isomorphism since $\dim_{\Bbbk} P_{S'} = \dim_{\Bbbk} M$. 

\item Suppose $\loewy(M) = n$ and $\pd(N) = 0 \text{ or } \infty $ if $\loewy(N) \leq n-1$. We can consider the following short exact sequence
$$0 \rightarrow S \rightarrow M \rightarrow N \rightarrow 0 $$
where $S = \soc(M)$, $\loewy(N) \leq n-1$ and $N$ does not have projective summands. By a similar argument to the previous part, we have a short exact sequence 
$$0 \rightarrow S \rightarrow P_{N} \rightarrow N \rightarrow 0$$
and an isomorphism $\pi: P_{N} \rightarrow M$.

\end{itemize}

\end{proof}

\end{lema}

\begin{coro}\label{5}
Consider $A = \frac{\Bbbk Q}{I}$ a bound quiver algebra. If $\phi_A(0)^{-1}$ is additive the subquiver $Q^{\infty}$ is closed by succesors.

\begin{proof}

The proof is straightforward from Lemmas \ref{1} and \ref{3}.

\end{proof}

\end{coro}

\begin{coro}\label{6}
Consider $A = \frac{\Bbbk Q}{I}$ a bound quiver algebra such that $\phi_A(0)^{-1}$ is an additive subcategory of $\mod A$. Then $P_{S_v}$ has simple socle for every vertex $v \in Q^{\infty}$.
\begin{proof}

Suppose there are at least two simple modules $S_1$ and $S_2$ such that $S_1 \oplus S_2 \subset \soc (P_{S_v})$. Since $Q^{\infty}$ is closed by succesors by Corollary \ref{5} then $S_1, S_2 \in \phi_A(0)^{-1}$. Consider the indecomposable modules $M_1 = \frac{P_{S_v}}{S_1}$, $M_2 = \frac{P_{S_v}}{S_2}$, $M_3 = \frac{P_{S_v}}{S_1 \oplus S_2}$. 
It is clear that $M_i \in \phi_A(0)^{-1}$ for $i = 1,2,3$. Since $\phi_A(0)^{-1}$ is additive $M= M_1\oplus M_2 \oplus M_3 \in \mod A$. 
On the other hand $M_1 \oplus M_2 \not \cong M_3$, and $\Omega (M_1 \oplus M_2) = \Omega(M_3)$. Then $\phi(M) \geq 1$, which is a contradiction. 

\end{proof}

\end{coro}

\begin{teo}\label{LIT maximal}
For a bound quiver algebra $A = \frac{\Bbbk Q}{I}$ are equivalent

\begin{enumerate}

\item $\phi_A(0)^{-1}$ is additive.
\item $\fidim \ \add\phi_A(0)^{-1} = 0$.
\item $\gl A < \infty$ or $A$ is selfinjective.

\end{enumerate}

\begin{proof}
\hfill 
\begin{itemize}
\item[($1 \Rightarrow 2$)] It is Lemma \ref{2}.
\item[($2 \Rightarrow 1$)] Suppose $M_1, M_2 \in \phi_A(0)^{-1}$. Since $\phi (M_1 \oplus M_2) = 0$, then $M_1\oplus M_2 \in \phi_A(0)^{-1}$.
\item[($1 \Rightarrow 3$)] By Corollary 11 of \cite{A}, $\gl A < \infty$ if and only if $\pd (S) < \infty$ for all $S \in \mathcal{S}_A$, so we can assume that there is a simple module $S$ with $\pd (S) = \infty$, then $Q^{\infty} \not = \emptyset$. 

If there are no simple modules with finite projective dimension, every module has either infinite projective dimension or it is a projective, $\phi_A(0)^{-1} = \mod A$. We conclude by Lemma \ref{4} that $A$ is a self-injective algebra. 

Suppose there are simple modules with finite projective dimension. Since $Q$ is a connected quiver and $Q^{\infty}$ is closed by sucessors, there is an arrow $\alpha$ with $\start(\alpha) = v_0$, $\target(\alpha) = w_0$ such that $\pd(S_{v_0}) < \infty$ and $\pd(S_{w_0}) = \infty$. We deduce that there is a vertex $w_1 \in Q^{\infty}$ such that $S_{w_1} \subset P_{S_{v_0}}$. On the other hand, by Corollary \ref{6} and the Pigeonhole principle there is a vertex $w_2 \in Q^{\infty}$ such that $S_{w_1} = \soc (P_{S_{w_2}})$. 
Now we can consider the indecomposable modules $M_1 = \frac{P_{S_{v_0}}}{S_{w_1}}$ and $M_2 = \frac{P_{S_{w_2}}}{S_{w_1}}$. It is clear that $M_1 \oplus M_2 \in \phi_{0}(A)$. 
On the other hand since $M_1 \not \cong M_2$ and $\Omega(M_1) = \Omega(M_2)$, $\phi(M_1 \oplus M_2) = 1$, which is contradiction. 

\item[($3 \Rightarrow 1$)] It is a particular case of Remark \ref{aditivas}.
\end{itemize}

\end{proof}

\end{teo}

As a straightforward consequence of the previous theorem, we have 

\begin{coro}
An algebra $A = \frac{\Bbbk Q}{I}$ is $(n,V,\phi^{-1}_A(0))$-LIT if and only if $\gl (A) < \infty$ or $A$ is selfinjective. As a consequence if $A$ is $(n,V,\phi^{-1}_A(0))$-LIT for some $n \in \mathbb{N}$ then $A$ is $(m, 0,\phi^{-1}_A(0))$-LIT for some $m \in \mathbb{N}$ 
\end{coro}

Since the global dimension and being selfinjective are left-right symmetric for Artin algebras. We have the following result

\begin{coro}\label{trivial} Let $A = \frac{\Bbbk Q}{I}$ be a bound quiver algebra. Then $\phi^{-1}_A(0)$ is 0-Igusa-Todorov if and only if $\phi^{-1}_{A^{op}}(0)$ is 0-Igusa-Todorov. As a consequence $A$ is $(n,0,\phi^{-1}_A(0))$-LIT for some $n \in \mathbb{N}$ if and only if $A^{op}$ is $(m,0,\phi^{-1}_{A^{op}}(0))$-LIT for some $m \in \mathbb{N}$.
\end{coro}


\section{Examples}

In this section, we see that the Corollary \ref{trivial} cannot be generalized to the family of LIT algebras, i.e. the notion of being a LIT algebra is not left-right symmetric.
The following example is based on Example 4 of \cite{BM}. For this reason, we omit part of the computation, the reader can find it in the cited article.

\begin{ej}\label{ejemplin}
Consider $\Bbbk$ an uncountable field. Let $A =\frac{\Bbbk Q}{I_A}$ be a finite dimensional $\Bbbk$-algebra and $B = \frac{\Bbbk Q'}{J_B^2}$, where $Q$ and $Q'$ are the following quivers
$$Q = \xymatrix{ & 0 \ar@(ul, dl)_{\gamma_1} \ar@(ur, dr)^{\gamma_2} \ar@(ru, lu)_{\gamma_3} },\ Q'=\xymatrix{ 1 \ar@(ul, dl)_{\bar{\beta}_1} \ar@/^2mm/[r]^{\beta_1} & 2 \ar@(ur, dr)^{\bar{\beta}_2} \ar@/^2mm/[l]^{\beta_2} }$$
and the ideal $I_A=\langle \gamma_i \gamma_i,\ \gamma_i \gamma_j + \gamma_j\gamma_i \ \forall \ i,j \in \{1,2,3\} \rangle$.\\

Consider $C = \frac{\Bbbk \Gamma}{I_C}$, with

\begin{itemize}

\item $\Gamma_0 = Q_0 \cup Q'_0$,

\item $\Gamma_1 = Q_1 \cup Q'_1 \cup \{ \alpha_0 : 0 \rightarrow 1\}$ and

\item $I_C = \langle I_A, J^2_{B} , \{ \lambda \alpha_0, \alpha_0\lambda \text{ } \forall \lambda \text{ such that }\length(\lambda)\geq 1 \} \rangle$. 

\end{itemize}

$C$ is not a LIT algebra (see \cite{BM}).\\

On the other hand, if we consider $C^{op}$. It is clear that $A \cong A^{op}$ and $B \cong B^{op}$ because its quivers and relations are symmetric respectively. So we conclude that $C^{op} = \frac{\Bbbk \Gamma}{I_{C^{op}}}$ with

\begin{itemize}

\item $\Gamma_0 = Q_0 \cup Q'_0$,

\item $\Gamma_1 = Q_1 \cup Q'_1 \cup \{ \alpha_0 : 1 \rightarrow 0\}$ and

\item $I_{C^{op}} = \langle I_A, J^2_{B} , \{ \lambda \alpha_0, \alpha_0\lambda \text{ } \forall \lambda \text{ such that }\length(\lambda)\geq 1 \} \rangle$. 

\end{itemize}

Since $A$ is selfinjective, then $\mod A$ is a $0$-Igusa-Todorov subcategory.
It is clear that $\Omega_C(\mod C) \subset \{ M \oplus S^{n_1}_1 \oplus S^{n_2}_2: M \in \mod A \subset \mod C, n_1, n_2 \in \mathbb{N}\}$. So we conclude that $C^{op}$ is $(1, \mod A, S_1\oplus S_2)$-LIT.

\end{ej}

\begin{obs}
The previous example also shows that
\begin{itemize}
\item in Proposition \ref{C LIT si A es IT y B LIT}, the hypothesis $A$ is an Igusa-Todorov cannot be changed for $A$ is a LIT algebra, and

\item in Theorem \ref{Algebra LIT triangular} (Vivero's theorem) the hypotesis on $M$ is necessary.

\end{itemize}
\end{obs}

The following remark shows that there exist $(n, V, \mathcal{D})$-LIT algebras where $\mathcal{P}_A \not \subset \mathcal{D} \not \subset \phi_{A}^{-1}(0)$.

\begin{obs}
Consider the algebra $C^{op}$ in Example \ref{ejemplin}.
\begin{itemize}
\item It is not an Igusa-Todorov algebra since $A$ is not an Igusa-Todorov algebra. 
\item Since $\gl (C^{op}) = \infty$ and $\Omega(S_2) = S_1 \oplus S_2 = \Omega(\frac{P_1}{S_1 \oplus S_2}) $, by Theorem \ref{LIT maximal}, the algebra $C^{op}$ is not a $(n, \phi_0^{-1}(A), \mathcal{D})$-LIT algebra for any $n \in \mathbb{N}$.
\item However we see in Example \ref{ejemplin} that $C^{op}$ is a $(1, \mod A, S_1\oplus S_2)$-LIT algebra. 

\end{itemize}
\end{obs}

\begin{obs} For LIT algebras, we have the following.

\begin{itemize}
\item Igusa-Todorov algebras $\not \subset \{A:$ $\phi_0^{-1}(A)$ is a 0-Igusa-Todorv subcategory $\}$.

Consider $A$ from Example \ref{ejemplin}. $A$ is a selfinjective algebra, then $A$ is an $(n,V,\phi_0^{-1}(A))$-LIT algebra. However $A$ is not an Igusa-Todorov algebra (see \cite{C}).

\item $\{A:$ $\phi_0^{-1}(A)$ is a 0-Igusa-Todorv subcategory $\} \not \subset$ Igusa-Todorov algebras.

On the other hand consider the algebra $C$ as in Proposition \ref{Igusa-Todorov cuadrada} where $A$ is the same as before, $Q_B$ is a vertex $v$, $\# J = 0$, $\# K = 1$ and $I_C = \langle I_A\rangle$. By Proposition \ref{Igusa-Todorov cuadrada} $C$ is Igusa-Todorov. Since $\Omega_C(S_v) = P_0$ and $\Omega_C(M) = \Omega_A(M)$ for every $A$-module $M$, then $\fidim (C) \geq 1$ and $\gl (C) = \infty$. Hence, by Theorem \ref{LIT maximal}, $C$ is not an $(n,V,\phi_0^{-1}(A))$-LIT algebra.

\end{itemize}

\end{obs}

\end{document}